\documentclass[a4paper,12pt]{article}
\usepackage{amsmath,amssymb,comment}

\newtheorem{prop}{Proposition}

\title{New examples of the representation of 1 by the sum of reciprocals of semiprime numbers}
\author{Tatsuru Watanabe \footnote{E-mail address: tat67wat@kje.biglobe.ne.jp}}
\date{\today}
\begin{document}

\maketitle
\begin{abstract}
 In $1978$, Allan.Wm.Johnson obtained an example of the representation of $1$ by the sum of reciprocals of
the product of two distinct prime numbers. His example has $48$ terms [1], and we had no examples 
which have less than $48$ terms until now. 
In this paper, we construct $17$ new examples that have less than $48$ terms. Since all of the new examples have $47$ terms and we have no examples which have less than $47$ terms, it is assumed that the minimum number of terms is $47$.  
\end{abstract} 

\section{Introduction}
In this note, we consider the representation of 1 by the finite sum of some reciprocals where each denominator is the product of two distinct primes.
These representations of 1 are equivalent to solutions of a Diophantine equation below:
\begin{align}
\sum_{i=1}^{n} \frac{1}{x_{i}}=1, \;\;\;\;  where \;\; 2\times 3 \leq x_{1} < x_{2}  < \dots < x_{n}
\end{align}
when $ \;\; x_{i} = p_{i} q_{i} $, $ \; p_{i} < q_{i} $ are primes for each $i$, $n$ is the number of terms.\\

\subsection{Previous research and our problems}
In 1977, Barbeau obtained a solution of (1) which has $101$ terms ($n=101$)[2].
In 1978, Allan.Wm.Johnson improved Barbeau's result and exhibited a solution of (1) which has $48$ terms ($n=48$)[1].\\
Richard K.Guy questioned whether the minimum number of terms is $48$ [3].
Here is Johnson's example:
\begin{align*}
	1&=\frac{1}{6}+\frac{1}{10}+\frac{1}{14}+\frac{1}{15}+\frac{1}{21}+\frac{1}{22}+\frac{1}{26}+\frac{1}{33}+\frac{1}{34}+\frac{1}{35}+\frac{1}{38}\\
&+\frac{1}{39}+\frac{1}{46}+\frac{1}{51}+\frac{1}{55}+\frac{1}{57}+\frac{1}{58}+\frac{1}{62}+\frac{1}{65}+\frac{1}{69}+\frac{1}{77}+\frac{1}{82}\\
&+\frac{1}{85}+\frac{1}{86}+\frac{1}{87}+\frac{1}{91}+\frac{1}{93}+\frac{1}{95}+\frac{1}{115}+\frac{1}{119}+\frac{1}{123}+\frac{1}{133}\\
&+\frac{1}{155}+\frac{1}{187}+\frac{1}{203}+\frac{1}{209}+\frac{1}{215}+\frac{1}{221}+\frac{1}{247}+\frac{1}{265}+\frac{1}{287}\\
&+\frac{1}{299}+\frac{1}{319}+\frac{1}{323}+\frac{1}{391}+\frac{1}{689}+\frac{1}{731}+\frac{1}{901}
\end{align*}
No examples which have $47$ or fewer terms are discovered.
In 2015, Nechemia Burshtein conjectured that there are other examples of (1) which have $48$ terms [4].

Our problem is the following two;
\begin{itemize}
\item To construct the new examples of (1) which has $47$ or less terms.($n \leq 47$)
\item To consider the minimum number of terms.
\end{itemize}
\section{Main results}
We have the new $17$ examples of (1) which have $47$ terms.
\section{Constructions}
\subsection{Reconsideration of Johnson's example}
We factorize all denominators of Johnson's example.
{\scriptsize
\begin{align*}
	\begin{aligned}
		6&= 2\times3,\\
		26&= 2\times13,\\
		46&= 2\times23,\\
		65&= 5\times13,\\
		87&= 3\times29,\\
		123&= 3\times41,\\
		215&= 5\times43,\\
		319&= 11\times29,\\
	\end{aligned}
	\quad
	\begin{aligned}
		10&= 2\times5,\\
		33&= 3\times11,\\
		51&= 3\times17,\\
		69&= 3\times23,\\
		91&= 7\times13,\\
		133&= 7\times19,\\
		221&= 13\times17,\\
		323&= 17\times19,\\
	\end{aligned}
	\quad
	\begin{aligned}
		14&= 2\times7,\\
		34&= 2\times17,\\
		55&= 5\times11,\\
		77&= 7\times11,\\
		93&= 3\times31,\\
		155&= 5\times31,\\
		247&= 13\times19,\\
		391&= 17\times23,\\
	\end{aligned}
	\quad
	\begin{aligned}
		15&= 3\times5,\\
		35&= 5\times7,\\
		57&= 3\times19,\\
		82&= 2\times41,\\
		95&= 5\times19,\\
		187&= 11\times17,\\
		265&= 5\times53,\\
		689&= 13\times53,\\
	\end{aligned}
        \quad
	\begin{aligned}
		21&= 3\times7,\\
		38&= 2\times19,\\
		58&= 2\times29,\\
		85&= 5\times17,\\
		115&= 5\times23,\\
		203&= 7\times29,\\
		287&= 7\times41,\\
		731&= 17\times43,\\
	\end{aligned}
        \quad
	\begin{aligned}
		22&= 2\times11,\\
		39&= 3\times13,\\
		62&= 2\times31,\\
		86&= 2\times43,\\
		119&= 7\times17,\\
		209&= 11\times19,\\
		299&= 13\times23,\\
		901&= 17\times53.
	\end{aligned}
\end{align*}
}
We call the product of two distinct primes as {\it a semiprime number}.\\
Since semiprime numbers exist infinitely, it is impossible to calculate all patterns under no restriction on denominators.\\
In the above table, all of the prime factors are less than $53$.
So in this study, we restrict that the maximum of the prime factors of denominators is $53$.

There are 120 semiprimes obtained by multiplying different two prime numbers from $2$ to $53$.

\subsection{Constructions of examples of $47$ terms}
We fix the vector ${\bf p}_{f}$ in $120$ dimensional vector space as follows:
$$
{\bf p}_{f}=( \frac{1}{6}, \frac{1}{10}, \frac{1}{14}, \frac{1}{15}, \dots , \frac{1}{2491}) \in \mathbb{R}^{120}
$$
This is a 120-dimensional vector having the reciprocals of semiprime numbers in descending order as components.
For ${\bf p}_{f} $, we correspond to the following another vector $ {\bf v} \in \mathbb {R} ^ {120} $:

$$
{\bf v} =( a_1, a_2, \dots ,a_{120} ) \in \mathbb{R}^{120}
$$
Each $ a_i $ takes a value of only $ 0 $ or $ 1 $. It takes $ 1 $ if the component of $ {\bf p} _ {f} $ appears in 
solutions of (1) , otherwise it takes $ 0 $.
For example, if $ \frac{1} {6} $ which is the first component of $ {\bf p} _ {f} $ appears in the representation, the first component $a_1$ of $ {\bf v} $ is equal to $ 1 $.
The number of reciprocals of semiprime numbers appearing in the solution of (1) is 
equal to the following inner product:
$$
A ={\bf v} \cdot { }^{t}(1,1, \dots ,1).
$$
Where $(1,1,\dots,1) \in \mathbb{R}^{120}$, all components of this vector are $1$.
Of course, we aim for A = $47$, but it is difficult, we decide in part the components of {\bf v}. 
Let ${\bf v} (n)$ be the vector determined up to the $n$-th component of ${\bf v}$ and all other components are $0$.
Consider the other inner product,
$$
B(n)= {\bf v} (n) \cdot { }^{t} {\bf p}_{f}.
$$
 $B(n)$ is a summation up to $n$-th terms.
We make 1 by adding $ B(n) $ and some reciprocals of semiprime numbers. 
We make a tree diagram for adding semiprime numbers or not, one by one to $B$, 
starting from the $n$-th smallest semiprime number according to the following proposition $1$.
 These observations are summarized in the following two propositions.
And we consider the fixed vector ${\bf q} _{f} \in \mathbb{R}^{120} $ where ${\bf q} _{f}=(6,10,14,15, \dots , 2491)$.
All of the components of  ${\bf q} _{f}$ are semiprime numbers which are the products of two different primes, its bigger prime factor is up to $53$. We assume that each component of ${\bf q} _{f}$ is arranged in ascending order. And each component of ${\bf q} _{f}$ is also reciprocal of corresponding component of ${\bf p} _{f}$.
\begin{prop}
We assume that $q_i$ represents the $i$-th component of  ${\bf q}_{f}$,\\
 $6=q_1, 10=q_2, 14=q_3, \dots 2491=q_{120} $.
\begin{description}
\item[1.1]{When the components of ${\bf v}$ are determined from $a_1$ to $a_n$,the inequality\\
 $ \displaystyle{ B(n)+ \sum_{k=1}^{47-C(n)} \frac{1}{q_{n+k}} < 1}$ does not hold, where $C(n)={\bf v}(n) \cdot  { }^{t}(1,1, \dots ,1)$}.
\item[1.2]{
For any prime numbers $r<53$, when the components of ${\bf v}$ are determined from $a_1$ to $a_n$,and the corresponding $n$-th component of the vector ${\bf q} _f$ is greater than $53r,$
multiples of $r$ do not appear in the subsequent semiprime numbers.
Therefore, the common denominator of the partial sum 
up to the above condition must not be a multiple of $r$.
}
\end{description}
\end{prop}

We found $2$ new examples of (1) consisting of $47$ terms, 
 and also we found $94$ examples consisting of $48$ terms with this restriction, including Johnson's example.

\subsection{Solutions of (1) which have $47$ terms}

Here are our new two examples of $47$ terms.
\begin{align*}
1&=\frac{1}{6}+\frac{1}{10}+\frac{1}{14}+\frac{1}{15}+\frac{1}{21}+\frac{1}{22}+\frac{1}{26}+\frac{1}{33}+\frac{1}{34}+\frac{1}{35}+\frac{1}{38}+\frac{1}{39}\\
&+\frac{1}{46}+\frac{1}{51}+\frac{1}{55}+\frac{1}{57}+\frac{1}{58}+\frac{1}{62}+\frac{1}{65}+\frac{1}{69}+\frac{1}{74}+\frac{1}{82}+\frac{1}{85}+\frac{1}{86}\\
&+\frac{1}{87}+\frac{1}{91}+\frac{1}{93}+\frac{1}{95}+\frac{1}{106}+\frac{1}{111}+\frac{1}{123}+\frac{1}{133}+\frac{1}{145}+\frac{1}{155}+\frac{1}{159}\\
&+\frac{1}{185}+\frac{1}{203}+\frac{1}{215}+\frac{1}{253}+\frac{1}{265}+\frac{1}{287}+\frac{1}{319}+\frac{1}{493}+\frac{1}{583}\\
&+\frac{1}{731}+\frac{1}{851}+\frac{1}{1073}\\
\end{align*}
\begin{align*}
1&=\frac{1}{6}+\frac{1}{10}+\frac{1}{14}+\frac{1}{15}+\frac{1}{21}+\frac{1}{22}+\frac{1}{26}+\frac{1}{33}+\frac{1}{34}+\frac{1}{35}+\frac{1}{38}+\frac{1}{39}\\
&+\frac{1}{46}+\frac{1}{51}+\frac{1}{55}+\frac{1}{57}+\frac{1}{58}+\frac{1}{62}+\frac{1}{65}+\frac{1}{74}+\frac{1}{77}+\frac{1}{82}+\frac{1}{85}+\frac{1}{87}\\
&+\frac{1}{91}+\frac{1}{93}+\frac{1}{95}+\frac{1}{111}+\frac{1}{115}+\frac{1}{119}+\frac{1}{123}+\frac{1}{129}+\frac{1}{133}+\frac{1}{143}\\
&+\frac{1}{145}+\frac{1}{155}+\frac{1}{161}+\frac{1}{221}+\frac{1}{253}+\frac{1}{259}+\frac{1}{287}+\frac{1}{299}+\frac{1}{391}\\
&+\frac{1}{473}+\frac{1}{481}+\frac{1}{1247}+\frac{1}{1591}
\end{align*}
Later, we shifted the upper bound of prime factors of denominators from $53$ to $101$.
The similar consideration and calculation in a $325$-dimensional vector space, with some minor revisions,
we obtained more $15$ examples consisting of $47$ terms. But there were no examples which have $46$ or fewer terms under this restriction.

\begin{align*}
1&=\frac{1}{6}+\frac{1}{10}+\frac{1}{14}+\frac{1}{15}+\frac{1}{21}+\frac{1}{22}+\frac{1}{26}+\frac{1}{33}+\frac{1}{34}+\frac{1}{35}+\frac{1}{38}+\frac{1}{39}\\
&+\frac{1}{46}+\frac{1}{51}+\frac{1}{55}+\frac{1}{57}+\frac{1}{58}+\frac{1}{62}+\frac{1}{65}+\frac{1}{69}+\frac{1}{74}+\frac{1}{77}+\frac{1}{86}\\
&+\frac{1}{87}+\frac{1}{91}+\frac{1}{93}+\frac{1}{94}+\frac{1}{95}+\frac{1}{111}+\frac{1}{115}+\frac{1}{118}+\frac{1}{119}+\frac{1}{129}+\frac{1}{133}\\
&+\frac{1}{145}+\frac{1}{155}+\frac{1}{161}+\frac{1}{185}+\frac{1}{187}+\frac{1}{329}+\frac{1}{517}+\frac{1}{667}+\frac{1}{851}+\frac{1}{1073}\\
&+\frac{1}{1357}+\frac{1}{1363}+\frac{1}{2537}\\
\end{align*}
\begin{align*}
1&=\frac{1}{6}+\frac{1}{10}+\frac{1}{14}+\frac{1}{15}+\frac{1}{21}+\frac{1}{22}+\frac{1}{26}+\frac{1}{33}+\frac{1}{34}+\frac{1}{35}+\frac{1}{38}+\frac{1}{39}\\
&+\frac{1}{46}+\frac{1}{51}+\frac{1}{55}+\frac{1}{57}+\frac{1}{58}+\frac{1}{62}+\frac{1}{65}+\frac{1}{69}+\frac{1}{77}+\frac{1}{82}+\frac{1}{85}\\
&+\frac{1}{86}+\frac{1}{87}+\frac{1}{91}+\frac{1}{93}+\frac{1}{94}+\frac{1}{95}+\frac{1}{118}+\frac{1}{119}+\frac{1}{123}+\frac{1}{133}+\frac{1}{145}\\
&+\frac{1}{155}+\frac{1}{161}+\frac{1}{187}+\frac{1}{203}+\frac{1}{215}+\frac{1}{287}+\frac{1}{329}+\frac{1}{493}+\frac{1}{517}+\frac{1}{1247}\\
&+\frac{1}{1357}+\frac{1}{1363}+\frac{1}{2537}\\
\end{align*}
\begin{align*}
1&=\frac{1}{6}+\frac{1}{10}+\frac{1}{14}+\frac{1}{15}+\frac{1}{21}+\frac{1}{22}+\frac{1}{26}+\frac{1}{33}+\frac{1}{34}+\frac{1}{35}+\frac{1}{38}+\frac{1}{39}\\
&+\frac{1}{46}+\frac{1}{51}+\frac{1}{55}+\frac{1}{57}+\frac{1}{58}+\frac{1}{62}+\frac{1}{65}+\frac{1}{69}+\frac{1}{77}+\frac{1}{82}+\frac{1}{85}\\
&+\frac{1}{86}+\frac{1}{87}+\frac{1}{93}+\frac{1}{94}+\frac{1}{95}+\frac{1}{115}+\frac{1}{118}+\frac{1}{119}+\frac{1}{123}+\frac{1}{129}+\frac{1}{133}\\
&+\frac{1}{141}+\frac{1}{143}+\frac{1}{145}+\frac{1}{155}+\frac{1}{287}+\frac{1}{377}+\frac{1}{391}+\frac{1}{551}+\frac{1}{799}+\frac{1}{893}\\
&+\frac{1}{1357}+\frac{1}{1363}+\frac{1}{2537}\\
\end{align*}
\begin{align*}
1&=\frac{1}{6}+\frac{1}{10}+\frac{1}{14}+\frac{1}{15}+\frac{1}{21}+\frac{1}{22}+\frac{1}{26}+\frac{1}{33}+\frac{1}{34}+\frac{1}{35}+\frac{1}{38}+\frac{1}{39}\\
&+\frac{1}{46}+\frac{1}{51}+\frac{1}{55}+\frac{1}{57}+\frac{1}{58}+\frac{1}{62}+\frac{1}{65}+\frac{1}{69}+\frac{1}{77}+\frac{1}{82}+\frac{1}{85}\\
&+\frac{1}{86}+\frac{1}{87}+\frac{1}{93}+\frac{1}{94}+\frac{1}{95}+\frac{1}{115}+\frac{1}{118}+\frac{1}{119}+\frac{1}{123}+\frac{1}{129}+\frac{1}{133}\\
&+\frac{1}{143}+\frac{1}{145}+\frac{1}{155}+\frac{1}{177}+\frac{1}{287}+\frac{1}{299}+\frac{1}{391}+\frac{1}{551}+\frac{1}{611}+\frac{1}{799}\\
&+\frac{1}{893}+\frac{1}{1711}+\frac{1}{2537}\\
\end{align*}
\begin{align*}
1&=\frac{1}{6}+\frac{1}{10}+\frac{1}{14}+\frac{1}{15}+\frac{1}{21}+\frac{1}{22}+\frac{1}{26}+\frac{1}{33}+\frac{1}{34}+\frac{1}{35}+\frac{1}{38}+\frac{1}{39}\\
&+\frac{1}{46}+\frac{1}{51}+\frac{1}{55}+\frac{1}{57}+\frac{1}{58}+\frac{1}{62}+\frac{1}{65}+\frac{1}{69}+\frac{1}{77}+\frac{1}{82}+\frac{1}{86}\\
&+\frac{1}{87}+\frac{1}{91}+\frac{1}{93}+\frac{1}{94}+\frac{1}{95}+\frac{1}{118}+\frac{1}{119}+\frac{1}{123}+\frac{1}{129}+\frac{1}{133}+\frac{1}{145}\\
&+\frac{1}{155}+\frac{1}{161}+\frac{1}{177}+\frac{1}{187}+\frac{1}{203}+\frac{1}{287}+\frac{1}{299}+\frac{1}{329}+\frac{1}{377}+\frac{1}{517}\\
&+\frac{1}{1363}+\frac{1}{1711}+\frac{1}{2537}\\
\end{align*}
\begin{align*}
1&=\frac{1}{6}+\frac{1}{10}+\frac{1}{14}+\frac{1}{15}+\frac{1}{21}+\frac{1}{22}+\frac{1}{26}+\frac{1}{33}+\frac{1}{34}+\frac{1}{35}+\frac{1}{38}+\frac{1}{39}\\
&+\frac{1}{46}+\frac{1}{51}+\frac{1}{55}+\frac{1}{57}+\frac{1}{58}+\frac{1}{62}+\frac{1}{65}+\frac{1}{69}+\frac{1}{74}+\frac{1}{77}+\frac{1}{82}\\
&+\frac{1}{86}+\frac{1}{91}+\frac{1}{93}+\frac{1}{95}+\frac{1}{111}+\frac{1}{119}+\frac{1}{122}+\frac{1}{123}+\frac{1}{133}+\frac{1}{143}+\frac{1}{145}\\
&+\frac{1}{155}+\frac{1}{161}+\frac{1}{183}+\frac{1}{187}+\frac{1}{259}+\frac{1}{287}+\frac{1}{299}+\frac{1}{319}+\frac{1}{473}+\frac{1}{481}\\
&+\frac{1}{559}+\frac{1}{671}+\frac{1}{1591}\\
\end{align*}
\begin{align*}
1&=\frac{1}{6}+\frac{1}{10}+\frac{1}{14}+\frac{1}{15}+\frac{1}{21}+\frac{1}{22}+\frac{1}{26}+\frac{1}{33}+\frac{1}{34}+\frac{1}{35}+\frac{1}{38}+\frac{1}{39}\\
&+\frac{1}{46}+\frac{1}{51}+\frac{1}{55}+\frac{1}{57}+\frac{1}{58}+\frac{1}{62}+\frac{1}{65}+\frac{1}{69}+\frac{1}{77}+\frac{1}{82}+\frac{1}{85}\\
&+\frac{1}{86}+\frac{1}{87}+\frac{1}{91}+\frac{1}{93}+\frac{1}{94}+\frac{1}{95}+\frac{1}{118}+\frac{1}{119}+\frac{1}{123}+\frac{1}{129}+\frac{1}{133}\\
&+\frac{1}{155}+\frac{1}{161}+\frac{1}{187}+\frac{1}{213}+\frac{1}{287}+\frac{1}{329}+\frac{1}{355}+\frac{1}{493}+\frac{1}{497}+\frac{1}{517}\\
&+\frac{1}{1357}+\frac{1}{1363}+\frac{1}{2537}\\
\end{align*}
\begin{align*}
1&=\frac{1}{6}+\frac{1}{10}+\frac{1}{14}+\frac{1}{15}+\frac{1}{21}+\frac{1}{22}+\frac{1}{26}+\frac{1}{33}+\frac{1}{34}+\frac{1}{35}+\frac{1}{38}+\frac{1}{39}\\
&+\frac{1}{46}+\frac{1}{51}+\frac{1}{55}+\frac{1}{57}+\frac{1}{58}+\frac{1}{62}+\frac{1}{65}+\frac{1}{69}+\frac{1}{74}+\frac{1}{77}+\frac{1}{82}\\
&+\frac{1}{85}+\frac{1}{87}+\frac{1}{91}+\frac{1}{93}+\frac{1}{95}+\frac{1}{115}+\frac{1}{119}+\frac{1}{123}+\frac{1}{129}+\frac{1}{133}+\frac{1}{143}\\
&+\frac{1}{145}+\frac{1}{155}+\frac{1}{187}+\frac{1}{213}+\frac{1}{221}+\frac{1}{287}+\frac{1}{301}+\frac{1}{559}+\frac{1}{629}+\frac{1}{781}\\
&+\frac{1}{1247}+\frac{1}{1591}+\frac{1}{1633}\\
\end{align*}
\begin{align*}
1&=\frac{1}{6}+\frac{1}{10}+\frac{1}{14}+\frac{1}{15}+\frac{1}{21}+\frac{1}{22}+\frac{1}{26}+\frac{1}{33}+\frac{1}{34}+\frac{1}{35}+\frac{1}{38}+\frac{1}{39}\\
&+\frac{1}{46}+\frac{1}{51}+\frac{1}{55}+\frac{1}{57}+\frac{1}{58}+\frac{1}{62}+\frac{1}{65}+\frac{1}{69}+\frac{1}{77}+\frac{1}{82}+\frac{1}{85}\\
&+\frac{1}{86}+\frac{1}{91}+\frac{1}{93}+\frac{1}{94}+\frac{1}{95}+\frac{1}{115}+\frac{1}{119}+\frac{1}{123}+\frac{1}{133}+\frac{1}{141}+\frac{1}{142}\\
&+\frac{1}{143}+\frac{1}{155}+\frac{1}{203}+\frac{1}{221}+\frac{1}{235}+\frac{1}{287}+\frac{1}{299}+\frac{1}{355}+\frac{1}{377}+\frac{1}{391}\\
&+\frac{1}{559}+\frac{1}{2021}+\frac{1}{2059}\\
\end{align*}
\begin{align*}
1&=\frac{1}{6}+\frac{1}{10}+\frac{1}{14}+\frac{1}{15}+\frac{1}{21}+\frac{1}{22}+\frac{1}{26}+\frac{1}{33}+\frac{1}{34}+\frac{1}{35}+\frac{1}{38}+\frac{1}{39}\\
&+\frac{1}{46}+\frac{1}{51}+\frac{1}{55}+\frac{1}{57}+\frac{1}{58}+\frac{1}{62}+\frac{1}{65}+\frac{1}{69}+\frac{1}{77}+\frac{1}{82}+\frac{1}{85}\\
&+\frac{1}{86}+\frac{1}{87}+\frac{1}{91}+\frac{1}{93}+\frac{1}{94}+\frac{1}{95}+\frac{1}{111}+\frac{1}{115}+\frac{1}{119}+\frac{1}{133}+\frac{1}{141}\\
&+\frac{1}{142}+\frac{1}{155}+\frac{1}{187}+\frac{1}{205}+\frac{1}{213}+\frac{1}{235}+\frac{1}{493}+\frac{1}{517}+\frac{1}{1363}+\frac{1}{1633}\\
&+\frac{1}{1763}+\frac{1}{1927}+\frac{1}{2627}\\
\end{align*}
\begin{align*}
1&=\frac{1}{6}+\frac{1}{10}+\frac{1}{14}+\frac{1}{15}+\frac{1}{21}+\frac{1}{22}+\frac{1}{26}+\frac{1}{33}+\frac{1}{34}+\frac{1}{35}+\frac{1}{38}+\frac{1}{39}\\
&+\frac{1}{46}+\frac{1}{51}+\frac{1}{55}+\frac{1}{57}+\frac{1}{58}+\frac{1}{62}+\frac{1}{65}+\frac{1}{69}+\frac{1}{77}+\frac{1}{82}+\frac{1}{86}\\
&+\frac{1}{91}+\frac{1}{93}+\frac{1}{94}+\frac{1}{95}+\frac{1}{111}+\frac{1}{115}+\frac{1}{119}+\frac{1}{123}+\frac{1}{133}+\frac{1}{141}+\frac{1}{142}\\
&+\frac{1}{143}+\frac{1}{145}+\frac{1}{155}+\frac{1}{203}+\frac{1}{213}+\frac{1}{221}+\frac{1}{235}+\frac{1}{287}+\frac{1}{493}+\frac{1}{559}\\
&+\frac{1}{1633}+\frac{1}{2021}+\frac{1}{2627}\\
\end{align*}
\begin{align*}
1&=\frac{1}{6}+\frac{1}{10}+\frac{1}{14}+\frac{1}{15}+\frac{1}{21}+\frac{1}{22}+\frac{1}{26}+\frac{1}{33}+\frac{1}{34}+\frac{1}{35}+\frac{1}{38}+\frac{1}{39}\\
&+\frac{1}{46}+\frac{1}{51}+\frac{1}{55}+\frac{1}{57}+\frac{1}{58}+\frac{1}{62}+\frac{1}{65}+\frac{1}{69}+\frac{1}{74}+\frac{1}{77}+\frac{1}{85}\\
&+\frac{1}{86}+\frac{1}{87}+\frac{1}{91}+\frac{1}{93}+\frac{1}{94}+\frac{1}{95}+\frac{1}{111}+\frac{1}{115}+\frac{1}{119}+\frac{1}{133}+\frac{1}{141}\\
&+\frac{1}{142}+\frac{1}{145}+\frac{1}{155}+\frac{1}{187}+\frac{1}{213}+\frac{1}{517}+\frac{1}{629}+\frac{1}{799}+\frac{1}{1247}+\frac{1}{1591}\\
&+\frac{1}{1633}+\frac{1}{2021}+\frac{1}{2627}\\
\end{align*}
\begin{align*}
1&=\frac{1}{6}+\frac{1}{10}+\frac{1}{14}+\frac{1}{15}+\frac{1}{21}+\frac{1}{22}+\frac{1}{26}+\frac{1}{33}+\frac{1}{34}+\frac{1}{35}+\frac{1}{38}+\frac{1}{39}\\
&+\frac{1}{46}+\frac{1}{51}+\frac{1}{55}+\frac{1}{57}+\frac{1}{58}+\frac{1}{62}+\frac{1}{65}+\frac{1}{69}+\frac{1}{77}+\frac{1}{82}+\frac{1}{85}\\
&+\frac{1}{87}+\frac{1}{91}+\frac{1}{93}+\frac{1}{94}+\frac{1}{95}+\frac{1}{111}+\frac{1}{115}+\frac{1}{118}+\frac{1}{119}+\frac{1}{123}+\frac{1}{133}\\
&+\frac{1}{142}+\frac{1}{155}+\frac{1}{187}+\frac{1}{213}+\frac{1}{287}+\frac{1}{295}+\frac{1}{329}+\frac{1}{413}+\frac{1}{493}+\frac{1}{517}\\
&+\frac{1}{1363}+\frac{1}{1633}+\frac{1}{2627}\\
\end{align*}
\begin{align*}
1&=\frac{1}{6}+\frac{1}{10}+\frac{1}{14}+\frac{1}{15}+\frac{1}{21}+\frac{1}{22}+\frac{1}{26}+\frac{1}{33}+\frac{1}{34}+\frac{1}{35}+\frac{1}{38}+\frac{1}{39}\\
&+\frac{1}{46}+\frac{1}{51}+\frac{1}{55}+\frac{1}{57}+\frac{1}{58}+\frac{1}{62}+\frac{1}{65}+\frac{1}{69}+\frac{1}{77}+\frac{1}{82}+\frac{1}{87}\\
&+\frac{1}{91}+\frac{1}{93}+\frac{1}{95}+\frac{1}{106}+\frac{1}{111}+\frac{1}{115}+\frac{1}{119}+\frac{1}{122}+\frac{1}{123}+\frac{1}{133}+\frac{1}{142}\\
&+\frac{1}{155}+\frac{1}{159}+\frac{1}{183}+\frac{1}{187}+\frac{1}{203}+\frac{1}{213}+\frac{1}{265}+\frac{1}{287}+\frac{1}{319}+\frac{1}{583}\\
&+\frac{1}{671}+\frac{1}{1633}+\frac{1}{2627}\\
\end{align*}
\begin{align*}
1&=\frac{1}{6}+\frac{1}{10}+\frac{1}{14}+\frac{1}{15}+\frac{1}{21}+\frac{1}{22}+\frac{1}{26}+\frac{1}{33}+\frac{1}{34}+\frac{1}{35}+\frac{1}{38}+\frac{1}{39}\\
&+\frac{1}{46}+\frac{1}{51}+\frac{1}{55}+\frac{1}{57}+\frac{1}{58}+\frac{1}{62}+\frac{1}{65}+\frac{1}{69}+\frac{1}{77}+\frac{1}{82}+\frac{1}{85}\\
&+\frac{1}{86}+\frac{1}{87}+\frac{1}{91}+\frac{1}{93}+\frac{1}{94}+\frac{1}{95}+\frac{1}{115}+\frac{1}{119}+\frac{1}{123}+\frac{1}{129}+\frac{1}{133}\\
&+\frac{1}{141}+\frac{1}{142}+\frac{1}{145}+\frac{1}{155}+\frac{1}{187}+\frac{1}{517}+\frac{1}{799}+\frac{1}{1207}+\frac{1}{1247}+\frac{1}{1633}\\
&+\frac{1}{1763}+\frac{1}{2021}+\frac{1}{2911}\\
\end{align*}
\section{Consideration of the minimum number of terms and future works}
\subsection{No $46$-term example}
We proved that there are no $46$-term examples when the maximum prime factor is $101$.
Also, the maximum prime factors of denominators of our $17$ examples which have $47$ terms, are $71$, so there is a low possibility of the existence of another example having $47$ terms except for our $17$ examples.
\subsection{Future works}
Our future works are to expand the dimension of vector space
and proof that the minimum number of terms is $47$.\\

Acknowledgments：
I would like to give heartwarming thanks to Kenta Yoshizaki who provided carefully considered feedback and valuable comments. 
I would also like to thank the Center for Computational Sciences at University of Tsukuba and all my teachers whose opinions and information have helped me very much throughout the production of this study.
And more than anyone, to my parents.


\begin{thebibliography}{99}

\bibitem[1]{[1]} Allan.WM.Johnson,{\it  Letter to editor}, Crux Mathematicorum (=Eureka (Ottawa)), 4 (1978) 190. 
\bibitem[2]{[2]} E.J.Barbeau, {\it Expressing one as a sum of distinct reciprocals }:
 Comments and a bibliography,Eureka (Ottawa) 3 (1977) p178-p181
\bibitem[3]{[3]} Richard K. Guy,{\it Unsolved Problems in Number Theory}
\bibitem[4]{[4]} Nechemia Burshtein, {\it On distinct unit fractions whose sum equals 1 when $x_i \mid \hspace{-.72em}/ x_j$ for $i \neq j$ }, Journal for 
Algebra and Number Theory Academia Volume 5, issue 4, Dec.2015, p117-p124.

\end{thebibliography}
\end{document}